# FRACTIONAL-ORDER STATE SPACE MODELS

**Ľubomír DORČÁK[1], Ivo PETRÁŠ[2], Imrich KOŠTIAL[3] and Ján TERPÁK[4]**
[1] Department of Informatics and Process Control
BERG Faculty, Technical University of Košice
B. Němcovej 3, 042 00 Košice, Slovak Republic
lubomir.dorcak@tuke.sk,
[2] ivo.petras@tuke.sk, [3] imrich.kostial@tuke.sk, [4] jan.terpak@tuke.sk

**Abstract:** In this paper we will present some alternative types of mathematical description and methods of solution of the fractional-order dynamical system in the state space. We point out the difference in the true sense of the name „state" space for the integer-order and fractional-order system and the importance of the initialization function for the fractional-order system. Some implications concerning the state feedback control theory are discussed. Presented are the results of simulations.

## 1 Introduction

Besides theoretical research of the fractional-order (FO) derivatives and integrals [Oldham, Spanier 1974, Samko, Kilbas, Marichev 1987, Podlubný 1994, 1999, Hartley, Lorenzo 1999, Podlubný 2000, Leszczynski, Ciesielski 2002 and many others] there are growing number of applications of the fractional calculus in such different areas as e.g. long electrical lines [Heaviside 1922], electrochemical processes, dielectric polarization [Westerlund 1994], colored noise [Mandelbrot 1967], viscoelastic materials, chaos, in control theory [Manabe 1961, Outstalup 1988, 1995, Axtel, Bise 1990, Kolojanov, Dimitrova 1992, Dorčák 1994, Podlubný, Dorčák, Koštial 1997] and in many other areas.

In works [Manabe 1961, Outstalup 1988, Axtel and Bise 1990, Kolojanov et al. 1992] the first generalizations of analysis methods for FO control systems were made (s-plane, frequency response, etc.). Some of our works were oriented to the methods of FO system parameters identification [Dorčák et al. 1996, Podlubný 1998], methods of FO controllers synthesis [Dorčák 1994, Petráš 1998, etc.], methods of stability analysis [Dorčák et al. 1998, Petráš et al. 1999], etc. Very interesting in control theory are the methods of solution of the vector FO differential equation or systems of FO differential equations, which is considered e.g. in works [Hartley et al. 1999, Dorčák et al. 2000, 2001, Leszczynski et al. 2002] and so on. Some new knowledge and problems appeared in this field.

So in this contribution we will present some alternative types of mathematical description and some alternative types of solution of the FO dynamical system in the state space. Compared are polynomial approximation and continued fraction expansion of Tustin's operator. We want to point out the difference in the true sense of the name „state" space for the integer order and FO system, the importance of the initialization function for the FO system, and some implications concerning the state feedback control theory.

## 2 Methods of description and solution of the FO systems and controllers

Fractional-order systems are described in time domain by a FO differential equation (DE) or by systems of FO differential equations (DEs). Very frequent in control theory, as a model of the controlled system, is the following differential equation

$$a_2 \, y^{(\alpha)}(t) + a_1 \, y^{(\beta)}(t) + a_0 \, y(t) = u(t) \qquad (1)$$

where $\alpha$, $\beta$ are generally real numbers, $a_2$, $a_1$, $a_0$ are arbitrary constants. The FO $PI^\lambda D^\delta$ controller is described by the FO DE $u(t) = Ke(t) + T_i \, D_t^{-\lambda} e(t) + T_d \, D_t^\delta e(t)$ where $\lambda$ and $\delta$ are integral and derivative orders, $K$, $T_i$, and $T_d$ are the proportional, integration, and derivation constant respectively.

Analytical solution of FO DEs [Podlubný 1994, 1999] is rather complicated and for higher order equations almost impossible. More convenient are numerical solutions. One is based on the Grunwald-Letnikov (GL) definition of the differintegral operator

$$_a D_t^r f(t) = \lim_{T \to 0} \frac{1}{T^r} \sum_{j=0}^{\left[\frac{t-a}{T}\right]} b_j^{(r)} f(t - jT) = \lim_{T \to 0} \frac{1}{T^r} \Delta_T^r \, f(t), \qquad (2)$$

where $[x]$ means the integer part $x$, $\Delta f(t)$ is the generalized finite difference of order $r$ and time step $T$, $b_j$ are binomial coefficients $b_0^{(r)} = 1, b_j^{(r)} = \left(1 - (1 + (\pm r))/j\right) b_{j-1}^{(r)}$. The discrete approximation of the fractional differintegral operator of order $r$ can be expressed by generating function $s = \omega(z^{-1})$, where $z^{-1}$ is the shift operator. For discrete approximation of time derivative, we can use the generating function corresponding to the Z-transform of backward difference rule $\omega(z^{-1}) = (1 - z^{-1})/T$. Performing the power series expansion (PSE) of $(1-z^{-1})^{\pm r}$ we obtain the Z version of the GL formula

$$D^{\pm r}(z) = (\omega(z^{-1}))^{\pm r} = T^{\mp r} z^{-[L/T]} \sum_{j=0}^{[L/T]} b_j^{(\pm r)} z^{[L/T]-j} \qquad (3)$$

where $L$ is the memory length. Using GL definition (2) or (3) we can derive the explicit equation for numerical solution of equation (1). In [Dorčák 1994] and other works this was done for many types of differential equations of control systems and these solutions were compared to analytical solutions of the equations with very satisfactory results.

If we want to describe an FO system, represented by equation (1), by a system of FO DEs, we ought to decompose equation (1). For this decomposition we should not neglect the general property of FO derivatives [Podlubný 1999], e.g. $({_aD_t^\alpha} y)(t) = D^m({_aI_t^{m-\alpha}} y)(t)$, ${_aI_t^\alpha}({_aI_t^\beta} y(t)) = {_aI_t^{\alpha+\beta}} y(t) = {_aI_t^{\beta+\alpha}} y(t)$ because, in general, the FO derivatives do not commute ${_aD_t^\alpha}({_aD_t^\beta} y(t)) \neq {_aD_t^{\alpha+\beta}} y(t) \neq {_aD_t^{\beta+\alpha}} y(t)$, except for zero initial conditions. The integer operator commutes with the fractional operator $(D_t^m({_aD_t^\alpha} y))(t) = ({_aD_t^{m+\alpha}} y)(t)$, but the reverse is truee. The mixed operators (derivatives and integrals) commute only in the following way $({_aD_t^\alpha}({_aI_t^\beta} y))(t) = ({_aD_t^{\alpha-\beta}} y)(t)$. A very interesting method of decomposition which uses the above properties is in [Leszczynski, Ciesielski 2002].

Almost all decomposition methods lead to the vector or state space representation expressing the FO derivatives on the left hand side of the system of equations. In our works [Dorčák et al. 2000, 2001] we were concerned with decomposition of some FO DEs under zero initial conditions expressing the integer order derivatives on the left hand side of the system of equations and FO derivatives on the right hand side. Now we will investigate the decomposition based on expressing the FO derivatives on the left hand side of the system of equations because of their better applicability in control theory.

Consider the system described by fractional-order differential equation (1) under zero initial conditions. After its modification and substitution of state space variables $x(t) = x_1(t)$, $x^{(\beta)}(t) = x_2(t)$ we can derive the following state space model

$$\begin{aligned} x_1^{(\beta)}(t) &= x_2(t) \\ x_2^{(\alpha-\beta)}(t) &= -\frac{a_0}{a_2} x_1(t) - \frac{a_1}{a_2} x_2(t) + \frac{1}{a_2} u(t) \\ y(t) &= x_1(t) \end{aligned} \qquad (4)$$

The structure of this system of the FO DEs is the same as for the integer order system (α=2, β=1) and the same as for the system of FO DEs obtained by expressing the FO derivatives on the right hand side of the system of equations [Dorčák et al. 2001]. Of course the results are the same too. After approximation of the FO derivatives in equations (4) by the GL formula (forward difference rule) we obtain equations for the discretized model

$$\begin{aligned} x_{1,k+1} &= -\sum_{j=1}^{k+1} b_j^{(\beta)} x_{1,k+1-j} + T^\beta x_{2,k} \\ x_{2,k+1} &= -\sum_{j=1}^{k+1} c_j^{(\alpha-\beta)} x_{2,k+1-j} + T^{\alpha-\beta}\left(-\frac{a_0}{a_2} x_{1,k} - \frac{a_1}{a_2} x_{2,k} + \frac{1}{a_2} u_k\right) \\ y_k &= x_{1,k} \end{aligned} \qquad (5)$$

As we see, FO DEs accumulate the whole information of the function in a weighted form. This is so called memory effect. Fractionally differentiated state space variables must be known as long as the system has been operated to obtain correct response. This is known as the initialization function [Hartley, Lorenzo 1999], for integer-order systems it is constant and for FO systems it is time varying. In the usual integer-order system theory, the set of states of the system, known at any given point in time, along with the system equations, are sufficient to predict the response of the system both forward or backward in time. The fractional dynamic variables do not represent the "state" of a system at any given time in the previous sense, we need all history of states, or sufficient number by short memory principle (minimum 100), for initialization function computation.

Many authors make decompositions to the system of FO DEs with the FO derivatives of order $r$ with limitation $0<r<1$. The consequence is that we have different systems of the FO DEs and different number of state space variables for the same equation (1) with only different α, β. But the results of these systems are equivalent. Therefore we use the above decomposition to system (4) without limitations to $r$ for different values of α, β.

Because of the above-mentioned memory effect and from this following high memory consumption we realized a more effective method of FO derivatives approximation based on direct discretization using the Tustin rule [Gorenflo 1996] and continued fraction expansion (CFE). The resulting discrete transfer function, approximating FO operators of order $r$ by the CFE of Tustin's rule, can be expressed [Vinagre et al. 2000] as

$$D^{\pm r}(z) = \frac{Y(z)}{U(z)} = \left(\frac{2}{T}\right)^{\pm r} \text{CFE}\left\{\left(\frac{1-z^{-1}}{1+z^{-1}}\right)^{\pm r}\right\}_{p,q} = \left(\frac{2}{T}\right)^{\pm r} \frac{P_p^r(z^{-1})}{Q_q^r(z^{-1})} \quad (6)$$

where $Y(z)$ and $U(z)$ are Z transforms of the output and input sequences $y(nT)$, $u(nT)$, $P$ and $Q$ are polynomials of degrees $p$ and $q$, respectively, in the variable $z^{-1}$. For $q=p=9$ we have $Q_9(z^{-1})$ = $(-52480r^3+147456r+r^9-120r^7+4368r^5)z^{-9}$ + $(45r^8+120330r^4-909765r^2-4410r^6+893025)z^{-8}$+$(-5742495r-76230r^5+1451835r^3+990r^7)z^{-7}$+$(-13097700+9514890r^2-796950r^4+13860r^6)z^{-6}$ + $(33648615r-5405400r^3+135135r^5)z^{-5}$ + $(-23648625r^2+51081030+945945r^4)z^{-4}$ + $(-61486425r+4729725r^3)z^{-3}$ + $(16216200r^2-72972900)z^{-2}$ + $(34459425r)z^{-1}$ +34459425. The polynomial $P_9$ has the same coefficients except the odd coefficients, which have opposite sign. Introducing this polynomials into equation (6), multiply this by $z$ (forward difference rule) and after substitution $R_i = z(2/T)^r Q_i$ we obtain from equations (4) the following equations

$$x_{1,k+1} = \frac{1}{R_0^\beta}\left(-\sum_{i=1}^{9} R_i^\beta x_{1,k+1-i} + \sum_{i=0}^{9} Q_i^\beta x_{2,k-i}\right)$$

$$x_{2,k+1} = \frac{1}{R_0^\delta}\left(-\frac{a_0}{a_2}\sum_{i=0}^{9} Q_i^\delta x_{1,k-i} - \frac{a_1}{a_2}\sum_{i=0}^{9} Q_i^\delta x_{2,k-i} - \sum_{i=1}^{9} R_i^\delta x_{2,k+1-i} + \frac{1}{a_2}\sum_{i=0}^{9} Q_i^\delta u_k\right) \quad (7)$$

$$y_k = x_{1,k}$$

where $\delta=\alpha-\beta$. The results of equations (7) are very close to results from equations (5) by at least ten times less memory requirements. In Fig. 1 and in Fig. 2 are depicted unit step re-sponses and corresponding state trajectories of equation (1) computed with PSE method (5) with L=100 and CFE method (7) with p=q=9, $a_2$=0.8, $a_1$=0.5, $a_0$=1.0, α=2.2, β=0.9, and T=0.1 [s] for both methods - only for comparison purposes. With equations (4) we have a possibility to describe model in vector and matrix relations. This system with above-mentioned parameters is controllable, because the rank of the matrix of controllability $Q_R$

$A = \begin{bmatrix} 0 & 1 \\ -1,25 & -0,625 \end{bmatrix}$, $B = \begin{bmatrix} 0 \\ 1,25 \end{bmatrix}$, $C = 1$, $Q_R \equiv Q_R^v = \begin{bmatrix} 0 & 1,25 \\ 1,25 & -0,7813 \end{bmatrix}$ is : rank($Q_R$) = 2. One method

for FO $PD^\delta$ and $PI^\lambda$ controller design was presented in work [Dorčák et al. 2001]. FO $PD^\delta$ controller gives us an extra degree of freedom, but it is dangerous, because we can obtain the controller with weak integrator and additional pole can cause system instability.

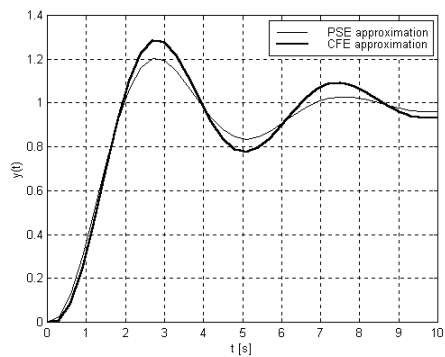 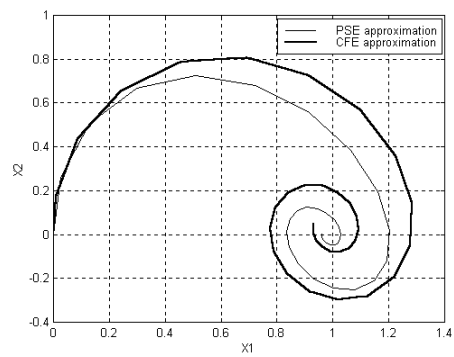

Figure 1. Unit step responses                Figure 2. State trajectories

## Acknowledgements

This work was partially supported by grant VEGA 1/7098/20 from the Slovak Grant Agency for Science.